\documentclass[aps,pra,twocolumn,showpacs,showkeys,superscriptaddress,nobalancelastpage,nofootinbib,floatfix,longbibliography]{revtex4-2}

\usepackage{amsfonts}
\usepackage{amsmath}
\usepackage{fancybox}
\usepackage{graphicx}
\usepackage{booktabs} 
\usepackage{siunitx}  
\usepackage{subcaption} 
\usepackage{tikz}
\usetikzlibrary{quantikz}
\usepackage{algorithm}
\usepackage{algpseudocode}
\usepackage{physics}
\usepackage{caption}
\usepackage{float}

\newcommand{\C}{\mathbb{C}}
\newcommand{\ii}{\mathrm{i}}

\newtheorem{remark}{Remark}[section]

\begin{document}

\title{Projection Coefficients Estimation in Continuous-Variable Quantum Circuits}

\author{M. W. AlMasri}
\email{mwalmasri2003@gmail.com}
\affiliation{Wilczek Quantum Center, School of Physics and Astronomy, Shanghai Jiao Tong University, Minhang, Shanghai, China}

\begin{abstract}
In this work, we propose a continuous-variable quantum algorithm to compute the projection coefficients of a holomorphic function in the Segal--Bargmann space by leveraging its isometric correspondence with single-mode quantum states. Using CV quantum circuits, we prepare the state $\ket{f}$ associated with $f(z)$ and extract the coefficients $c_n = \braket{n}{f}$ via photon-number-resolved detection, enhanced by interferometric phase referencing to recover full complex amplitudes. We detail the construction of the state-preparation oracle for various functional classes and analyze the protocol's robustness under realistic noise models, including detector inefficiency and state preparation errors. This enables direct quantum estimation and visualization of the coefficient sequence---offering a hardware-native protocol for characterizing non-Gaussian states and analyzing functions defined by quantum oracles, complementary to classical numerical integration.
\end{abstract}

\maketitle


\section{Introduction}
Series expansions provide a foundational tool for approximating and analyzing functions in mathematical physics, enabling the reduction of nonlinear or transcendental problems to tractable polynomial or spectral forms. Central to this framework is the Maclaurin series—the Taylor expansion about $z_{0}=0$—which, for a function $f : \mathbb{C} \to \mathbb{C}$ holomorphic in a neighborhood of $z = 0$, takes the form  \cite{Ahlfors,Fokas}
\begin{equation}
f(z) = \sum_{n=0}^{\infty} \frac{f^{(n)}(0)}{n!}\, z^n,
\end{equation}
where $f^{(n)}(0) = \left. \frac{d^n f}{dz^n} \right|_{z=0}$ denotes the $n$-th complex derivative at the origin. The series converges absolutely and uniformly on every compact subset of the disk of analyticity $\{ z \in \mathbb{C} : |z| < R \}$, where $R > 0$ is the radius of convergence, determined by the distance to the nearest singularity in the complex plane. For multivariate functions $f : \mathbb{C}^d \to \mathbb{C}$, the Maclaurin expansion generalizes via multi-index notation:
\begin{eqnarray}
f(\mathbf{z}) = \sum_{\alpha \in \mathbb{N}_0^d}
\frac{1}{\alpha!} \,
\partial^\alpha f(\mathbf{0})\,
\mathbf{z}^\alpha,
\\
\alpha! = \alpha_1! \cdots \alpha_d!, \quad
\mathbf{z}^\alpha = z_1^{\alpha_1} \cdots z_d^{\alpha_d},
\end{eqnarray}
with $\partial^\alpha = \partial^{|\alpha|}/(\partial z_1^{\alpha_1} \cdots \partial z_d^{\alpha_d})$ and $|\alpha| = \alpha_1 + \cdots + \alpha_d$. In physical applications such as perturbative quantum mechanics, effective field theories, or semiclassical approximations, the $n$-th order truncation
\begin{equation}
f(\mathbf{z}) \approx \sum_{|\alpha| \le n}
\frac{1}{\alpha!} \,
\partial^\alpha f(\mathbf{0})\,
\mathbf{z}^\alpha
\end{equation}
yields controlled approximations whose error is governed by the remainder term (e.g., Lagrange or integral form). Crucially, for holomorphic functions, the Maclaurin series is not merely asymptotic: it converges to $f$ identically within its domain of analyticity, reflecting the deep link between complex differentiability and infinite differentiability \cite{Fokas}. \\

In this work, we introduce the projection coefficients algorithm in the Segal--Bargmann space $\mathcal{H}_{SB}$, a space of holomorphic functions equipped with a Gaussian integration measure, to find the leading-order series expansion terms graphically \cite{Segal,Bargmann,Hall}. The Segal--Bargmann space $\mathcal{H}_{SB}$, or Bargmann representation has been studied in a variety of physical systems   \cite{Perelomov,Stenholm,Bishop,almasri,master,Wahiddin, neuron,complex, central}. The accuracy of the computed derivative values depends on the precision of the numerical routine used to approximate the inner products. Furthermore, we present a continuous-variable (CV) quantum circuit implementation of the Projection Coefficients Algorithm, leveraging the isomorphism between holomorphic functions in the Segal--Bargmann space and single-mode quantum states \cite{PRL,review,Menicucci,Gu,Adesso}. By preparing the state $\ket{f}$ corresponding to a target function $f(z)$ and performing photon-number-resolved measurements---augmented with interferometric phase referencing---we directly extract the projection coefficients $c_n = \braket{n}{f}$ as measurable quantum amplitudes, enabling both estimation and visualization of their magnitude.

\section{Analytic Hilbert Spaces with Gaussian Measure}
In an analytic Hilbert space $\mathcal{H}$, the inner product between two complex functions $ f(z) $ and $ g(z) $ is typically given by:
\begin{equation}
\langle f, g \rangle = \int_{\Omega} \overline{f(z)} g(z) w(z) \, d\mu(z),
\end{equation}
where, $\Omega$ is  the domain of analyticity (e.g., the complex plane $\mathbb{C}$, the unit disk $\mathbb{D}$, etc.), $w(z)$ is the  weight function that ensures convergence of the integral and defines the specific Hilbert space, and  $d\mu(z)$ is the measure on the domain $\Omega$ (e.g., area measure $dA(z)$, arc-length measure $|dz|$, etc.).  The Segal--Bargmann space $\mathcal{H}_{SB}$  is a Hilbert space of entire functions $ f(z) $ that are square-integrable with respect to the Gaussian measure:
\begin{equation}
\|f(z)\|^2 = \int_{\mathbb{C}} |f(z)|^2 e^{-|z|^2} \, dA(z),
\end{equation}
where $ dA(z) $ is the area measure on the complex plane $ \mathbb{C} $. The monomials $ \{z^n\}_{n=0}^\infty $ form an orthogonal basis for $ B $, with norms given by:
\begin{equation}
\|z^n\|^2 = \pi n!.
\end{equation}
In the Segal--Bargmann space $\mathcal{H}_{SB}$ of entire functions, the inner product is \cite{Bargmann,Segal}:
\begin{equation}\label{innerb}
\langle f, g \rangle = \int_{\mathbb{C}} \overline{f(z)} g(z) e^{-|z|^2} \, dA(z),
\end{equation}
where, $\Omega = \mathbb{C}$ is the complex plane, $w(z) = e^{-|z|^2}$ is the Gaussian weight, and $d\mu(z) = dA(z)$ is the area measure. To ensure convergence of the integral the following conditions must be satisfied:
\begin{enumerate}
\item \textit{Holomorphicity of $f(z)$}: $f(z)$ must be an entire function (holomorphic everywhere in $\mathbb{C}$).
\item \textit{Growth condition on $f(z)$}: $f(z)$ must satisfy a growth bound such as:
\begin{equation}
|f(z)| \leq C e^{A|z|^2}, \quad \text{for some constants } C > 0 \text{ and } A < 1.
\end{equation}
This ensures that $f(z)e^{-|z|^2}$ decays exponentially as $|z| \to \infty$.
The polynomial factor does not affect convergence because the Gaussian factor $e^{-|z|^2}$ dominates.
\end{enumerate}
If these conditions are satisfied, the integral converges.\vskip 2mm

Let $\mathcal{F} = \text{span}\{\,|n\rangle\,\}_{n=0}^\infty$ be the Fock space of a single bosonic mode, where $|n\rangle$  is the number states (Fock states), $a$ is the annihilation operator,  and  $ a^\dagger$ is the creation operator. The actions of the annihilation and creation operators $a$ and $a^\dagger$ are given by:
\begin{align}
a|n\rangle &= \sqrt{n} |n - 1\rangle, \\
a^\dagger|n\rangle &= \sqrt{n + 1} |n + 1\rangle.
\end{align}
The Bargmann map $\Psi: \mathcal{F} \to \mathcal{H}_{SB}$ is defined by its action on the basis vectors as:
\begin{equation}
\Psi(|n\rangle) = z^n.
\end{equation}
Thus, any state $|\psi\rangle = \sum_{n=0}^\infty c_n |n\rangle$ is mapped to a function:
\begin{equation}
f(z) = \sum_{n=0}^\infty c_n z^n,
\end{equation}
which is an entire function on $\mathbb{C}$, and belongs to the Segal--Bargmann space $\mathcal{H}_{SB}$ equipped with the inner product Eq. \ref{innerb}.
Under this mapping, the annihilation operator $a$ corresponds to the derivative $\frac{\partial}{\partial z}$, the creation operator $a^\dagger$ corresponds to multiplication by $z$, and the number operator $a^\dagger a$ becomes $z \frac{\partial}{\partial z}$. \vskip 5mm

In the context of spin systems, Bargmann  representation can be constructed using the Schwinger-Jordan oscillator model , which maps spin operators onto bosonic creation and annihilation operators \cite{almasri}.
The Schwinger-Jordan mapping expresses spin operators $\mathbf{S} = (S_x, S_y, S_z)$ in terms of bosonic creation ($a^\dagger, b^\dagger$) and annihilation ($a, b$) operators. Using the transformation relations:
\begin{equation}
a \mapsto \frac{\partial}{\partial z_a}, \quad b \mapsto \frac{\partial}{\partial z_b}.
\end{equation}
and inner product of the form:
\begin{equation}
\langle f | g \rangle = \int_{\mathbb{C}^2} \overline{f(z_a, z_b)} g(z_a, z_b) e^{-|z_a|^2 - |z_b|^2} d\mu(z_a, z_b),
\end{equation}
where, $d\mu(z_a, z_b) = \frac{1}{\pi^2} dx_a dy_a dx_b dy_b$ is the normalized measure in the complex plane, and  $z_a = x_a + i y_a$ and $z_b = x_b + i y_b$.
Thus, the Bargmann representation encodes the spin states of the system as holomorphic functions $f(z_a, z_b)$ defined on the complex plane. Below, we discuss spin coherent states and Fock basis states as prototypical examples.  \vskip 2mm

A spin coherent state $|\alpha\rangle$ is defined as:
\begin{equation}
|\alpha\rangle = e^{\alpha S_- - \alpha^* S_+} |s, s\rangle,
\end{equation}
where $|s, s\rangle$ is the highest-weight state with $S_z = s$. In the Bargmann representation, this becomes:
\begin{equation}
f_{\text{coh}}(z_a, z_b) = \left( z_a + \alpha z_b \right)^{2s},
\end{equation}
where $z_a$ and $z_b$ are the Bargmann variables corresponding to the two bosonic modes.
\vskip 2mm

The Fock basis states $|n_a, n_b\rangle$ correspond to monomials in the Bargmann representation:
\begin{equation}
f_{n_a, n_b}(z_a, z_b) = \frac{z_a^{n_a} z_b^{n_b}}{\sqrt{n_a! n_b!}},
\end{equation}
with the constraint $n_a + n_b = 2s$.
This inner product ensures that the Bargmann representation is unitarily equivalent to the original Fock space representation.
\vskip 5mm

The Schwinger-Jordan mapping expresses spin operators $\mathbf{S} = (S_x, S_y, S_z)$ in terms of bosonic creation ($a^\dagger, b^\dagger$) and annihilation ($a, b$) operators. For a spin system with total spin $s$, the spin operators $S_+, S_-, S_z$ take simple forms in the Bargmann representation \cite{almasri}:
\begin{equation}
S_+ = a^\dagger b \mapsto z_a \frac{\partial}{\partial z_b}.
\end{equation}
\begin{equation}
S_- = b^\dagger a \mapsto z_b \frac{\partial}{\partial z_a}.
\end{equation}
\begin{equation}
S_z = \frac{1}{2}(a^\dagger a - b^\dagger b) \mapsto \frac{1}{2}\left( z_a \frac{\partial}{\partial z_a} - z_b \frac{\partial}{\partial z_b} \right).
\end{equation}
These operators act directly on the holomorphic functions $f(z_a, z_b)$.

\section{ Projection Coefficients }
The monomials $\{z^n\}_{n=0}^\infty$ constitute an orthogonal basis of the Segal--Bargmann space $\mathcal{H}_\mathrm{SB}$, endowed with the Gaussian inner product. Consequently, any $f \in \mathcal{H}_\mathrm{SB}$ admits a unique expansion
\begin{equation}
f(z) = \sum_{n=0}^\infty c_n z^n,
\qquad
c_n = \frac{\langle z^n, f \rangle}{\|z^n\|^2}.
\end{equation}
Using the Maclaurin series $f(z) = \sum_{k=0}^\infty a_k z^k$ and orthogonality $\langle z^n, z^k \rangle = \pi n!\, \delta_{nk}$, one finds
\begin{equation}
\langle z^n, f \rangle = a_n \pi n!, \qquad \|z^n\|^2 = \pi n!,
\end{equation}
hence
\begin{equation}\label{eq:cn_an}
c_n = a_n.
\end{equation}
Thus, the Segal--Bargmann projection coefficients coincide exactly with the Maclaurin coefficients—mirroring how Fourier coefficients arise from $L^2$-projections in harmonic analysis. The Gaussian weight $e^{-|z|^2}$ ensures convergence and induces a reproducing kernel structure, rendering $\mathcal{H}_\mathrm{SB}$ a Hilbert space of entire functions.
Below, we propose a general algorithm for computing the projection coefficients $c_{n}$ and test it using holomorphic test functions and truncation error analysis.

\subsection{Projection Coefficients Algorithm}
\textbf{Input:}
- A holomorphic function $ f(z) $.
- Maximum degree $ N $ (up to which to compute $ c_n $).
- Parameters for numerical integration:
- $ \texttt{num\_points} $: Number of grid points for Riemann sum approximation.
- $ \texttt{radius} $: Radius of the disk in the complex plane for integration.

\textbf{Output:}
- List of projection coefficients $ c_0, c_1, \dots, c_N $.

\textbf{Steps:}
\begin{enumerate}
\item Initialize an empty list to store the projection coefficients: $ \texttt{coefficients} = [\; ] $.
\item For each $ n = 0, 1, \dots, N $:
\begin{enumerate}
\item Compute the norm-squared of $ z^n $ in the Segal--Bargmann space $\mathcal{H}_{\text{SB}}$:
$$
\|z^n\|^2 = \pi n!
$$
\item Generate a polar grid in the complex plane:
- Radial coordinate: $ r \in [0, \texttt{radius}] $, discretized into $ \texttt{num\_points} $.
- Angular coordinate: $ \theta \in [0, 2\pi] $, discretized into $ \texttt{num\_points} $.
\item Approximate the inner product using a Riemann sum:
$$
\langle z^n, f(z) \rangle \approx \sum_{r, \theta} \overline{\left( r e^{i\theta} \right)^n} f(r e^{i\theta}) e^{-r^2} r \, \Delta r \, \Delta \theta,
$$
where $ \Delta r = \frac{\texttt{radius}}{\texttt{num\_points}} $ and $ \Delta \theta = \frac{2\pi}{\texttt{num\_points}} $.
\item Alternatively, when applicable and higher accuracy is desired, approximate the inner product using double quadrature (\texttt{dblquad}) to compute:
$$
\langle z^n, f(z) \rangle = \int_0^\infty \int_0^{2\pi} \overline{\left( r e^{i\theta} \right)^n} f\left(r e^{i\theta}\right) e^{-r^2} r \, d\theta \, dr.
$$
\item Compute the projection coefficient:
$$
c_n = \frac{\langle z^n, f(z) \rangle}{\|z^n\|^2}.
$$
\item Append $ c_n $ to the list of coefficients.
\end{enumerate}
\item Return the list of projection coefficients $ c_0, c_1, \dots, c_N $.
\end{enumerate}

\begin{figure}[htbp]
\centering
\includegraphics[width=8cm]{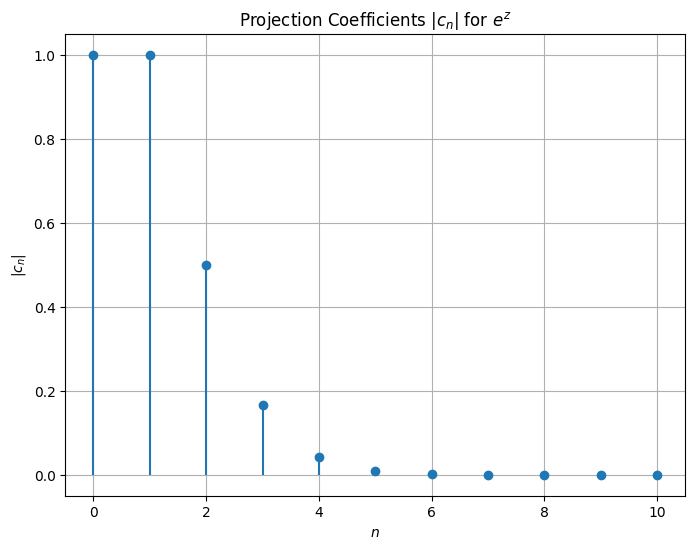}
\caption{The projection coefficients $c_{n}$ with $n\in \{0,10\}$ for the function $f(z)=e^{z}$. }
\label{fig:ez}
\end{figure}

\begin{table}[h!]
\centering
\begin{tabular}{cccc}
\toprule
$ n $ & $ \mathrm{Re}(c_n) $ & $ \mathrm{Im}(c_n) $ & $ |c_n| $ \\
\midrule
0 & 0.99999989 & 0.00000000 & 0.99999989 \\
1 & 0.99999809 & 0.00000000 & 0.99999809 \\
2 & 0.49999184 & 0.00000000 & 0.49999184 \\
3 & 0.16665114 & 0.00000000 & 0.16665114 \\
4 & 0.04164998 & 0.00000000 & 0.04164998 \\
5 & 0.00832180 & 0.00000000 & 0.00832180 \\
6 & 0.00138332 & 0.00000000 & 0.00138332 \\
7 & 0.00019643 & 0.00000000 & 0.00019643 \\
8 & 0.00002426 & 0.00000000 & 0.00002426 \\
9 & 0.00000264 & 0.00000000 & 0.00000264 \\
10 & 0.00000025 & 0.00000000 & 0.00000025 \\
\bottomrule
\end{tabular}
\caption{Numerical results of the projection coefficients $ c_n $ for $ e^z $.}
\label{projection_coefficients}
\end{table}

For illustration purposes, Let us consider the elementary  functions such as the exponential function $ f(z) = e^z $. The projection coefficients are plotted in Figure~\ref{fig:ez}. To compute the inner products, we used the double quadrature method (\texttt{dblquad}) from the \texttt{SciPy} library \cite{scipy}. This approach demonstrated significantly better accuracy compared to the Riemann sum. The numerical values of the projection coefficients are provided in Table~\ref{projection_coefficients}. \\

\begin{figure}[htbp]
\centering
\includegraphics[width=0.8\linewidth]{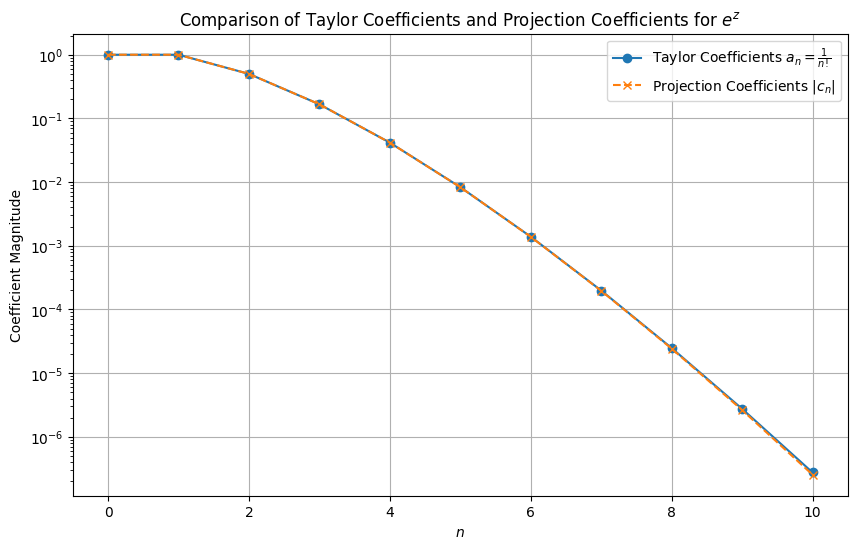}
\caption{Comparison of Maclaurin Coefficients and Projection Coefficients for $e^{z}$ and $n\in \{0,10\}$. }
\label{comparsion}
\end{figure}

The Maclaurin series expansion of $ e^z $ is given by:
\begin{equation}
e^z = \sum_{n=0}^\infty a_n z^n, \quad \text{where } a_n = \frac{1}{n!}.
\end{equation}
The Maclaurin coefficients for this function are listed in Table~\ref{taylor_coefficients} for $ n \in \{0, 10\} $. We plot the values from both table~\ref{projection_coefficients} and table~\ref{taylor_coefficients}. It is evident from Figure \ref{comparsion}that the Taylor coefficients and the projection coefficients are identical with remarkable accuracy.

\begin{table}[h!]
\centering
\begin{tabular}{cc}
\toprule
$ n $ & $ a_n = \frac{1}{n!} $ \\
\midrule
0 & 1.00000000 \\
1 & 1.00000000 \\
2 & 0.50000000 \\
3 & 0.16666667 \\
4 & 0.04166667 \\
5 & 0.00833333 \\
6 & 0.00138889 \\
7 & 0.00019841 \\
8 & 0.00002480 \\
9 & 0.00000276 \\
10 & 0.00000028 \\
\bottomrule
\end{tabular}
\caption{Maclaurin coefficients $ a_n = \frac{1}{n!} $ for the expansion of $ e^z $.}
\label{taylor_coefficients}
\end{table}

\begin{figure}[htbp]
\centering
\includegraphics[width=0.8\linewidth]{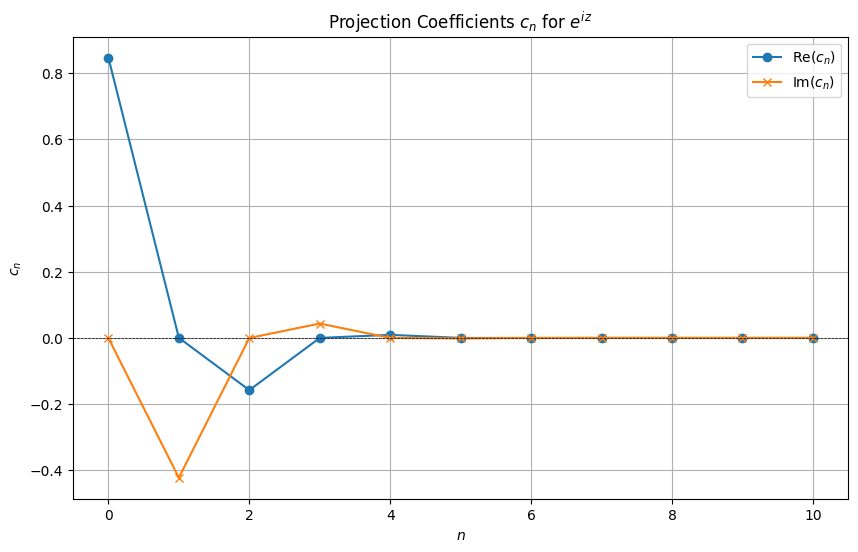}
\caption{The imaginary and real parts of projection coefficients for $f(z)=e^{iz}$.}
\label{eiz1}
\end{figure}

Similarly , consider the holomorphic function $f(z)=e^{iz}$ to be the test function. We plot the  projection coefficients $c_{n}$ with $n\in \{0,10\}$ in \ref{eiz1}. The Maclaurin series expansion of $ e^{iz} $ is given by:
$$
e^{iz} = \sum_{n=0}^\infty a_n z^n, \quad \text{where } a_n = \frac{(i)^n}{n!}.
$$

\begin{table}[h!]
\centering
\begin{tabular}{cccc}
\hline
$ n $ & $ \mathrm{Re}(c_n) $ & $ \mathrm{Im}(c_n) $ & $ |c_n| $ \\
\hline
0   &  0.84666501       & -0.00040490      & 0.84666511 \\
1   &  0.00027459       & -0.42244138      & 0.42244147 \\
2   & -0.15781570       & -0.00016986      & 0.15781579 \\
3   &  0.00001754       &  0.04339763      & 0.04339763 \\
4   &  0.00928797       & -0.00002342      & 0.00928800 \\
5   & -0.00000055       & -0.00159853      & 0.00159853 \\
6   & -0.00022354       & -0.00000190      & 0.00022355 \\
7   & -0.00000017       &  0.00002544      & 0.00002545 \\
8   &  0.00000247       & -0.00000011      & 0.00000247 \\
9   & -0.00000001       & -0.00000023      & 0.00000023 \\
10  & -0.00000002       & -0.00000000      & 0.00000002 \\
\hline
\end{tabular}
\caption{Numerical results of the projection coefficients $ c_n $ for $ f(z) = e^{iz} $.}
\label{projection_eiz}
\end{table}

\begin{table}[h!]
\centering
\begin{tabular}{cccc}
\hline
$ n $ & $ i^n $ & $ \mathrm{Re}(a_n) $ & $ \mathrm{Im}(a_n) $ \\
\hline
0 & $ 1 $     &  1.00000000 &  0.00000000 \\
1 & $ i $     &  0.00000000 &  1.00000000 \\
2 & $ -1 $    & -0.50000000 &  0.00000000 \\
3 & $ -i $    &  0.00000000 & -0.16666667 \\
4 & $ 1 $     &  0.04166667 &  0.00000000 \\
5 & $ i $     &  0.00000000 &  0.00833333 \\
6 & $ -1 $    & -0.00138889 &  0.00000000 \\
7 & $ -i $    &  0.00000000 & -0.00019841 \\
8 & $ 1 $     &  0.00002480 &  0.00000000 \\
9 & $ i $     &  0.00000000 &  0.00000276 \\
10 & $ -1 $   & -0.00000028 &  0.00000000 \\
\hline
\end{tabular}
\caption{Maclaurin coefficients $ a_n = \frac{(i)^n}{n!} $ for $ e^{iz} $.}
\label{taylor_coefficients_eiz}
\end{table}

The projection coefficients $ c_n $ are computed for the holomorphic function $ f(z) = e^{iz} $. Table~\ref{projection_eiz} displays the real part ($ \text{Re}(c_n) $), imaginary part ($ \text{Im}(c_n) $), and magnitude ($ |c_n| $) of these coefficients. The exact values are provided in table~\ref{taylor_coefficients_eiz}. By comparing table~\ref{projection_eiz} with table~\ref{taylor_coefficients_eiz}, we find the results to be highly informative regarding the order of the leading terms in the Maclaurin  series expansion. We employed the Riemann sum approximation to compute the integrals over the Segal--Bargmann space. However, it is reasonable to consider alternative numerical methods that offer improved accuracy and precision. \vskip 5mm

\begin{figure}[htbp]
\centering
\begin{subfigure}[b]{0.48\textwidth}
\includegraphics[width=0.8\linewidth]{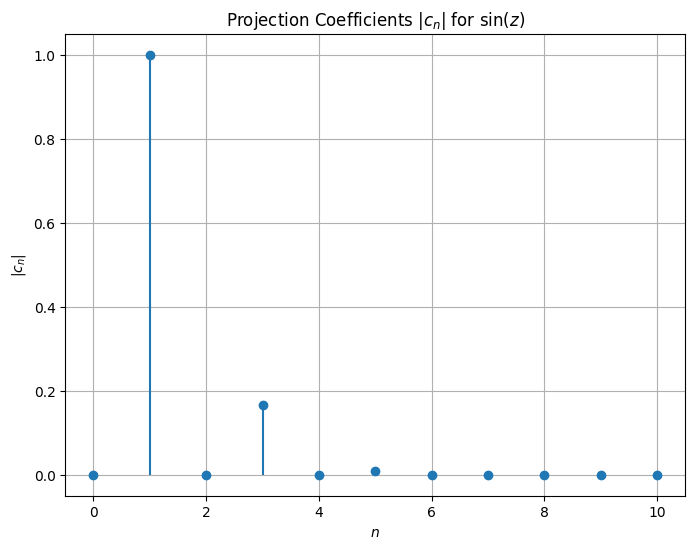}
\caption{}
\label{fig:sinz}
\end{subfigure}
\hfill
\begin{subfigure}[b]{0.48\textwidth}
\includegraphics[width=0.8\linewidth]{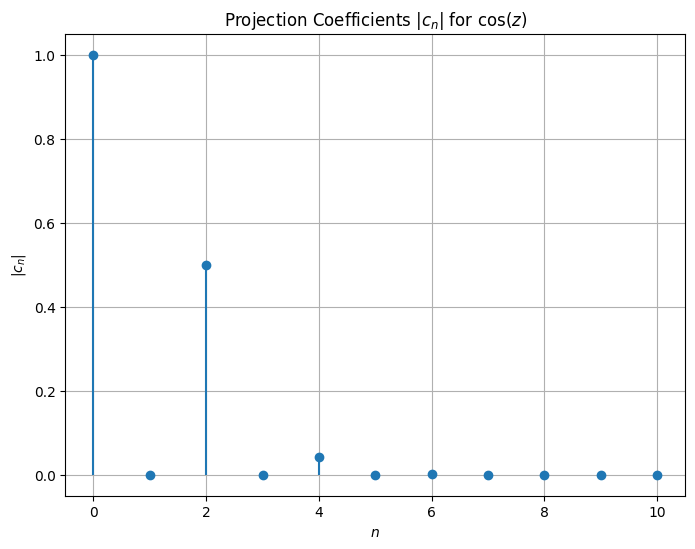}
\caption{}
\label{fig:cosz}
\end{subfigure}
\caption{The magnitudes of the projection coefficients are shown for (a) $\sin(z)$ and (b) $\cos(z)$. For $\sin(z)$, the leading terms occur at $n=1$ and $n=3$. For $\cos(z)$, the leading terms appear at $n=0$ and $n=2$. }
\label{fig:sin_cos_comparison}
\end{figure}

One interesting feature of the projection coefficients graphs is that they allow one to easily infer certain properties of functions, such as whether they are odd or even, as well as identify the leading terms in their Maclaurin series expansions. This enables safe and accurate approximations. To illustrate this, we consider two elementary functions one odd and one even as test cases.
In Figure~\ref{fig:sin_cos_comparison}, we plot the projection coefficients for $\sin(z)$ and $\cos(z)$. It is evident that for $\sin(z)$, the magnitudes $|c_0| = 0$ and $|c_2| = 0$, indicating that $\sin(z)$ is an odd function. Conversely, for $\cos(z)$, we observe that $|c_1| = 0$ and $|c_3| = 0$, confirming that $\cos(z)$ is an even function.
Furthermore, Figure~\ref{fig:sin_cos_comparison} suggests that $\sin(z)$ can be effectively approximated by retaining  the terms corresponding to $n=1,3$ in its Maclaurin series expansion. Similarly, $\cos(z)$ can be approximated by keeping the terms with $n=0,2,4$ in its Maclaurin series expansion.

To verify the conclusions drawn from Figure~\ref{fig:sin_cos_comparison}, we expand the functions $\sin(z)$ and $\cos(z)$ using their Maclaurin series. The Maclaurin series expansions of $\sin(z)$ and $\cos(z)$ around $z = 0$ are given by:
\begin{equation}
\sin(z) = \sum_{n=0}^\infty \frac{(-1)^n}{(2n+1)!} z^{2n+1} = z - \frac{z^3}{3!} + \frac{z^5}{5!} - \frac{z^7}{7!} + \cdots,
\end{equation}
and
\begin{equation}
\cos(z) = \sum_{n=0}^\infty \frac{(-1)^n}{(2n)!} z^{2n} = 1 - \frac{z^2}{2!} + \frac{z^4}{4!} - \frac{z^6}{6!} + \cdots.
\end{equation}
For practical purposes, these functions can be approximated by truncating the series after a few terms. More specifically:\\ 
- The approximation of $\sin(z)$ up to the third-order term is:
\begin{equation}
\sin(z) \approx z - \frac{z^3}{6} = c_1 z +c_3 z^3,
\end{equation}
where $c_1 = 1$ and $c_3 = -\frac{1}{6}$.\\ 
- The approximation of $\cos(z)$ up to the fourth-order term is:
\begin{equation}
\cos(z) \approx 1 - \frac{z^2}{2} + \frac{z^4}{24} = c_0 +c_2 z^2 + c_4 z^4,
\end{equation}
where $c_0 = 1$, $c_2 = -\frac{1}{2}$, and $c_4 = \frac{1}{24}$.
These approximations are accurate for small values of $z$, as higher-order terms become negligible when $|z|$ is small.

\subsection{Truncation Error Analysis}
\begin{figure}[htbp]
\centering
\includegraphics[width=0.9\linewidth]{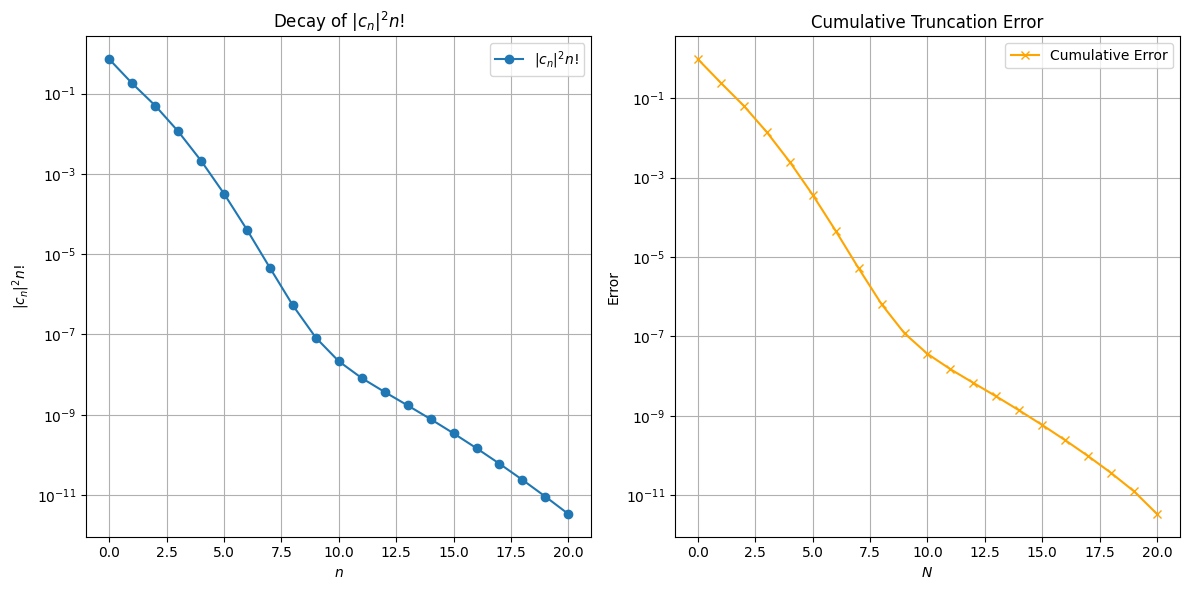}
\caption{The decay of $ |c_n|^2 n!$ and cumulative truncation error  for the function $e^{z}$ with $n\in \{0,20\}$. }
\label{truncation}
\end{figure}

If we approximate $ f(z) $ by a finite sum of terms up to degree $ N $, we have:
\begin{equation}
f_N(z) = \sum_{n=0}^N c_n z^n,
\end{equation}
The truncation error is the difference between $ f(z) $ and its approximation $ f_N(z) $. Specifically, the truncation error $ E_N(z) $ is:
\begin{equation}
E_N(z) = f(z) - f_N(z).
\end{equation}
In terms of the projection coefficients, the truncation error can be expressed as:
\begin{equation}
E_N(z) = \sum_{n=N+1}^\infty c_n z^n.
\end{equation}
To quantify the truncation error, we compute its norm in the Segal--Bargmann space $\mathcal{H}_{SB}$. The squared norm of $ E_N(z) $ is:
\begin{equation}
\|E_N(z)\|^2 = \int_{\mathbb{C}} |E_N(z)|^2 e^{-|z|^2} \, dA(z).
\end{equation}
Substituting $ E_N(z) = \sum_{n=N+1}^\infty c_n z^n $, we have:
\begin{equation}
|E_N(z)|^2 = \left| \sum_{n=N+1}^\infty c_n z^n \right|^2.
\end{equation}
Using the orthogonality of the monomials $ z^n $ in the Segal--Bargmann space $\mathcal{H}_{SB}$, the norm simplifies to:
\begin{equation}
\|E_N(z)\|^2 = \sum_{n=N+1}^\infty |c_n|^2 \|z^n\|^2.
\end{equation}
Since $ \|z^n\|^2 = \pi n! $, this becomes:
\begin{equation}
\|E_N(z)\|^2 = \pi \sum_{n=N+1}^\infty |c_n|^2 n!.
\end{equation}
The truncation error depends on the magnitude of the projection coefficients $ c_n $ and the factorial growth of $ n! $. If $ |c_n|^2 $ decays sufficiently fast as $ n \to \infty $, the truncation error will become small for large $ N $. Conversely, if $ f(z) $ has significant contributions from high-order terms (large $ n $), the truncation error may remain substantial even for moderately large $ N $. In figure \ref{truncation}, we plot the $ |c_n|^2 n!$ and cumulative truncation error  for the function $e^{z}$.

\section{CV Quantum Circuit Implementation}
The Segal--Bargmann (Bargmann--Fock) space provides an isometric isomorphism between holomorphic functions  and pure states of a single bosonic mode:
\begin{eqnarray}
f(z) \;\longleftrightarrow\; \ket{f} \coloneqq \sum_{n=0}^{\infty} c_n \ket{n},
\\
c_n = \braket{n}{f} = \frac{1}{\sqrt{\pi n!}} \int_{\C} \overline{z}^n f(z) e^{-|z|^2} \dd^2 z.
\label{eq:bargmann_isometry}
\end{eqnarray}
Therefore, the projection coefficients $\{c_n\}$ are precisely the Fock-basis amplitudes of the quantum state $\ket{f}$. The goal of the Projection Coefficients Algorithm is to estimate $\{c_n\}_{n=0}^N$ for a target $f(z)$ and visualize their magnitude/phase—e.g., to analyze non-Gaussianity, convergence, or analytic structure.
We propose a fully quantum implementation that prepares $\ket{f}$ (or an approximation thereof) on a CV platform and extracts $\{c_n\}$ via measurement.

\subsection{State Preparation Oracle}
Assume access to a unitary state-preparation oracle $\hat{U}_f$ such that
\begin{equation}
\hat{U}_f \ket{0} = \ket{f} = \frac{1}{\sqrt{\mathcal{N}}} \sum_{n=0}^{\infty} c_n \ket{n},
\qquad \mathcal{N} = \sum_{n=0}^{\infty} |c_n|^2 < \infty.
\label{eq:Uf}
\end{equation}
For holomorphic $f$, $\hat{U}_f$ may be constructed variationally (e.g., using a layered CV ansatz~\cite{killoran}), or via engineered photon addition for polynomial targets.

\subsubsection{Oracle Construction and Circuit Complexity}
The construction of the state-preparation oracle $\hat{U}_f$ is contingent upon the analytic structure of the target function $f(z)$. We identify three distinct regimes, each requiring specific quantum resources and offering different trade-offs between circuit depth and experimental feasibility.

\textbf{(i) Polynomial targets (Finite Truncation):} 
For functions represented by a finite polynomial $f(z) = \sum_{n=0}^N c_n z^n$, the corresponding quantum state $\ket{f}$ resides in a finite-dimensional subspace of the Fock space spanned by $\{\ket{0}, \dots, \ket{N}\}$. Such states can be prepared deterministically (up to heralding probability) using sequential photon-addition operations \cite{Zavatta2004}. Starting from the vacuum $\ket{0}$, we apply a cascade of non-Gaussian operators:
\begin{equation}
\hat{U}_f = \prod_{n=1}^{N} \left( \cos\theta_n \, \mathbb{I} + \sin\theta_n \, \hat{a}^\dagger \right),
\end{equation}
where the mixing angles $\theta_n$ are calibrated to reproduce the target coefficients $c_n$.  This approach requires exactly $N$ non-Gaussian steps, resulting in a circuit depth scaling linearly as $O(N)$. While conceptually straightforward, the experimental realization relies on high-efficiency heralded photon addition (e.g., via spontaneous parametric down-conversion), which can be resource-intensive due to probabilistic generation. However, for low-order truncations ($N \leq 5$), this method offers exact state preparation without variational error.

\textbf{(ii) Exponential-type functions (Coherent States):} 
A special but highly relevant class of functions takes the form $f(z) = e^{\alpha z}$, which maps directly to the coherent state $\ket{\alpha}$ in the Bargmann representation (ignoring normalization constants absorbed in the measure). For these targets, the oracle simplifies to the displacement operator:
\begin{equation}
\hat{U}_f = \hat{D}(\alpha) = \exp\left(\alpha \hat{a}^\dagger - \alpha^* \hat{a}\right).
\end{equation}
This is a purely Gaussian operation realizable with a single optical element (e.g., interfering the vacuum with a strong classical local oscillator on a highly transmissive beam splitter). \\
\textbf{(iii) General holomorphic functions (Variational Ansatz):} 
For arbitrary entire functions (e.g., $\cos z$, $\sin z$, or functions without closed-form expansions) that require infinite series truncated at $N$, deterministic construction is often intractable. Instead, we employ a variational continuous-variable quantum circuit ansatz \cite{killoran}:
\begin{equation}
\hat{U}_f(\vec{\phi}) = \prod_{l=1}^{L} \hat{R}(\phi_{l,1}) \hat{S}(\phi_{l,2}) \hat{D}(\phi_{l,3}),
\end{equation}
where $\hat{R}$, $\hat{S}$, and $\hat{D}$ denote rotation, squeezing, and displacement gates, respectively. The parameter vector $\vec{\phi} \in \mathbb{R}^{3L}$ is optimized via a classical feedback loop to minimize the infidelity loss function:
\begin{equation}
\mathcal{L}(\vec{\phi}) = 1 - \left| \braket{f_{\text{target}}}{\psi(\vec{\phi})} \right|^2,
\label{eq:infidelity_loss}
\end{equation}
where $\ket{\psi(\vec{\phi})} = \hat{U}_f(\vec{\phi})\ket{0}$. The optimization proceeds using one of two complementary strategies, selected based on the noise profile and differentiability of the hardware landscape:

\begin{enumerate}
    \item \textbf{Gradient Descent with Parameter-Shift Rule:} 
    For Gaussian gates, analytic gradients can be estimated without finite-difference approximations using the parameter-shift rule . For a parameter $\phi_k$ associated with a generator $\hat{G}_k$ (e.g., $\hat{n}$ for rotation, $\hat{x}\hat{p} + \hat{p}\hat{x}$ for squeezing), the gradient component is:
    \begin{equation}
    \frac{\partial \mathcal{L}}{\partial \phi_k} = \frac{1}{2} \left[ \mathcal{L}\left(\vec{\phi} + \frac{\pi}{2}\hat{e}_k\right) - \mathcal{L}\left(\vec{\phi} - \frac{\pi}{2}\hat{e}_k\right) \right],
    \end{equation}
    where $\hat{e}_k$ is the unit vector in the $k$-th direction. The parameters are updated via $\vec{\phi}^{(t+1)} = \vec{\phi}^{(t)} - \eta \nabla \mathcal{L}(\vec{\phi}^{(t)})$, with learning rate $\eta$. This method offers quadratic convergence near minima but requires $2 \times 3L$ circuit evaluations per step and is sensitive to shot noise unless $M \gg 1/\eta^2$ shots are used per evaluation.
    
    \item \textbf{Nelder-Mead Simplex Method:} 
    In regimes where gradient estimation is prohibitive due to high shot noise or non-differentiable hardware constraints, we employ the derivative-free Nelder-Mead algorithm \cite{nelder1965simplex}. The optimizer maintains a simplex of $3L+1$ vertices $\{\vec{\phi}^{(0)}, \dots, \vec{\phi}^{(3L)}\}$ in parameter space. At each iteration, the vertex with highest loss $\vec{\phi}_h$ is replaced via geometric operations:
    \begin{itemize}
        \item \textit{Reflection:} $\vec{\phi}_r = \vec{\phi}_c + \alpha(\vec{\phi}_c - \vec{\phi}_h)$,
        \item \textit{Expansion:} If $\mathcal{L}(\vec{\phi}_r) < \mathcal{L}(\vec{\phi}_l)$, try $\vec{\phi}_e = \vec{\phi}_c + \gamma(\vec{\phi}_r - \vec{\phi}_c)$,
        \item \textit{Contraction:} If reflection fails, contract towards the centroid $\vec{\phi}_c$,
        \item \textit{Shrink:} If all else fails, shrink all vertices towards the best point $\vec{\phi}_l$,
    \end{itemize}
    with standard coefficients $\alpha=1, \gamma=2, \rho=0.5, \sigma=0.5$. This method is robust against stochastic noise in $\mathcal{L}$ but typically converges slower than gradient-based approaches ($\mathcal{O}(N^2)$ vs.\ $\mathcal{O}(N)$ evaluations per step).
\end{enumerate}

\begin{remark}[Convergence and Initialization]
To ensure global convergence, we initialize $\vec{\phi}^{(0)}$ near the identity operation ($\vec{\phi} \approx 0$) with small random perturbations to break symmetry. The optimization terminates when either:
\begin{equation}
\mathcal{L}(\vec{\phi}^{(t)}) < \epsilon_{\text{tol}} \quad \text{or} \quad \|\vec{\phi}^{(t+1)} - \vec{\phi}^{(t)}\|_2 < \delta_{\text{step}},
\end{equation}
with typical thresholds $\epsilon_{\text{tol}} = 10^{-3}$ and $\delta_{\text{step}} = 10^{-4}$. For noisy intermediate-scale quantum (NISQ) devices, we recommend adaptive shot allocation, increasing $M$ as $\mathcal{L}$ approaches $\epsilon_{\text{tol}}$ to reduce variance in the final gradient estimates.
\end{remark}
In our numerical simulations for states truncated at $N \leq 10$, we observe that a circuit depth of $L \approx 3N$ layers typically suffices to achieve fidelities $\mathcal{F} > 0.95$. This approach leverages universal Gaussian gates (which are high-fidelity and scalable in time-domain multiplexed architectures \cite{Menicucci2014}) combined with a minimal number of non-Gaussian resources if higher-order non-classicality is required. The primary bottleneck shifts from gate count to the classical optimization convergence time, which remains manageable for low-photon-number subspaces.\vskip 5mm

\textbf{Experimental feasibility:} Current CV platforms based on time-domain multiplexed optical architectures \cite{Menicucci2014} have demonstrated scalability to large mode numbers and high-fidelity Gaussian operations in experimental implementations. For the benchmark functions considered ($e^z$, $\cos z$), the required circuit depth ($L \leq 30$) is well within near-term experimental capabilities. The primary bottleneck is PNR detector efficiency, addressed in the following noise analysis.

\subsection{Quantum Measurement of $c_n$}
Since $c_n = \braket{n}{f}$, the coefficients can be obtained by Fock basis tomography \cite{tomography,tomography1}. The simplest protocol uses repeated projective measurements in the number basis:
\begin{enumerate}
\item Prepare $\ket{f}$ via $\hat{U}_f \ket{0}$.
\item Perform photon-number-resolving detection (PNRD) on the mode \footnote{PNRD is a measurement that projects the state onto the Fock basis, yielding outcome $n$ with probability $p_n = |\braket{n}{\psi}|^2$.In real experiments, PNRD devices (e.g., transition-edge sensors, multiplexed SNSPDs) are non-ideal:
finite efficiency $\eta < 1$ $\rightarrow$ detected distribution is a smeared version of $p_n$, limited dynamic range (e.g., reliable only up to $n \sim 10$)\cite{Hadfield}.
}
\item Record outcome $n$.
\item Repeat $M$ times to estimate the probability distribution
\begin{equation}
p_n = |\braket{n}{f}|^2 = \frac{|c_n|^2}{\mathcal{N}}.
\label{eq:pn}
\end{equation}
\end{enumerate}
The magnitudes are then $|c_n| = \sqrt{\mathcal{N} \, p_n}$. To recover the phases $\arg(c_n)$, we use a phase-reference interferometer (Fig.~\ref{fig:circuit}):

\begin{figure}[htbp]
\centering
\begin{quantikz}[row sep=0.3cm, column sep=0.5cm]
\lstick[wires=2]{$\ket{0,0}$}
& \gate[2]{\hat{U}_{\mathrm{BS}}(\pi/4)}
& \gate{\hat{U}_f} & \qw & \meter{$n$} \\
& & \gate{\hat{D}(\alpha)} & \qw & \meter{$k$}
\end{quantikz}
\caption{Interferometric circuit for phase-sensitive estimation of $c_n$. A coherent state $\ket{\alpha}$ is prepared on the ancilla (via $\hat{D}(\alpha)\ket{0}$) and interfered with $\ket{f}$ on a 50:50 beam splitter $\hat{U}_{\mathrm{BS}}(\pi/4) = \exp\!\big[\tfrac{\pi}{4}(\hat{a}\hat{b}^\dagger - \hat{a}^\dagger \hat{b})\big]$. Joint PNRD yields statistics sensitive to $\arg(c_n)$.}
\label{fig:circuit}
\end{figure}

The joint probability of detecting $(n,k)$ photons is
\begin{equation}
P_{n,k}(\alpha) = \big| \bra{n,k} \hat{U}_{\mathrm{BS}}(\tfrac{\pi}{4}) \big( \ket{f} \otimes \ket{\alpha} \big) \big|^2.
\end{equation}
Expanding $\ket{\alpha} = e^{-|\alpha|^2/2} \sum_{m=0}^\infty \frac{\alpha^m}{\sqrt{m!}} \ket{m}$ and using the beam-splitter transformation, one finds (to leading order in $1/\sqrt{M}$):
\begin{equation}
\sqrt{P_{n,k}(\alpha)} \propto \left| c_n \frac{\alpha^k}{\sqrt{k!}} + c_{n+k} \frac{(-\alpha^*)^n}{\sqrt{n!}} + \cdots \right|.
\end{equation}
By fixing $k=0$ and scanning $\alpha = r e^{\ii\theta}$, the interference term $c_n \alpha^0 + c_{n+0}(-\alpha^*)^n$ isolates $\arg(c_n)$ via sinusoidal modulation in $\theta$.
In practice, for moderate $N$ ($\leq 8$), phase can be reconstructed via maximum-likelihood estimation over $\{P_{n,k}(\alpha_j)\}_{j=1}^J$ for a set of displacements $\{\alpha_j\}$.

\subsection{Algorithm Summary}
The full quantum Projection Coefficients Algorithm proceeds as:
\begin{algorithm}[H]
\caption{CV Quantum Projection Coefficients Algorithm }
\begin{algorithmic}[1]
\State \textbf{Input:} Oracle $\hat{U}_f$, max degree $N$, shots $M$, displacements $\{\alpha_j\}_{j=1}^J$.
\For{$n = 0$ to $N$}
\State Prepare $\ket{f} = \hat{U}_f \ket{0}$
\State Measure in Fock basis $M$ times $\Rightarrow$ estimate $p_n$
\State Compute $|c_n| \gets \sqrt{p_n}$ (assuming normalized $\ket{f}$)
\EndFor
\For{$j = 1$ to $J$}
\State Prepare $\ket{f} \otimes \ket{\alpha_j}$
\State Interfere on 50:50 BS
\State Record joint PNRD outcomes $\{(n,k)\}$
\EndFor
\State Reconstruct $\arg(c_n)$ via phase-fitting or ML
\State \textbf{Output:} $\{c_n = |c_n| e^{\ii \arg(c_n)}\}_{n=0}^N$
\State Plot:  $|c_n|$ vs $n$
\end{algorithmic}
\end{algorithm}

\subsection{Visualization}
The coefficients are typically visualized as:
\begin{itemize}
\item A bar plot of $|c_n|^2$ (occupation probabilities),
\item A polar plot of $c_n$ in the complex plane (revealing analytic structure),
\item A semi-log plot of $|c_n|$ to assess decay (e.g., Gaussian $\sim 1/\sqrt{n!}$, exponential $\sim r^n$).
\end{itemize}
For example, for $f(z) = e^{\alpha z - \frac{1}{2}|\alpha|^2}$ (coherent state), one obtains $c_n = \alpha^n / \sqrt{n!}$, yielding a Poisson-weighted rotation in the complex plane—verified experimentally via the above circuit.
This implementation leverages the intrinsic link between the Bargmann representation and Fock space, turning an analytical projection problem into a hardware-native quantum measurement protocol suitable for low-to-moderate photon numbers.

To validate the quantum measurement protocol under realistic conditions, we simulate the estimation of projection coefficient magnitudes $|c_n|$ for $f(z) = \cos(z)$ incorporating experimental noise sources. The exact coefficients are analytically known: $c_n = (-1)^{n/2}/n!$ for even $n$, and $c_n = 0$ for odd $n$. We construct the normalized quantum state $\ket{f} \propto \sum_{n=0}^{N} c_n \ket{n}$ (with truncation at $N=10$), then simulate projective Fock measurements with the following noise model:

\begin{enumerate}
\item \textit{Detector inefficiency:} Modeled as a beam splitter with transmission $\eta$ preceding ideal PNR detection \cite{Hadfield}. The detected photon number $n_{\text{det}}$ follows a binomial distribution $P(n_{\text{det}}|n) = \binom{n}{n_{\text{det}}} \eta^{n_{\text{det}}} (1-\eta)^{n-n_{\text{det}}}$.
\item \textit{Dark counts:} Added Poissonian noise with mean $\bar{n}_{\text{dark}} = 10^{-3}$ per measurement window.
\item \textit{State preparation infidelity:} Gaussian-distributed amplitude errors with standard deviation $\sigma_{\text{prep}} = 0.02$.
\end{enumerate}

\begin{figure}[htbp]
\centering
\includegraphics[width=0.75\linewidth]{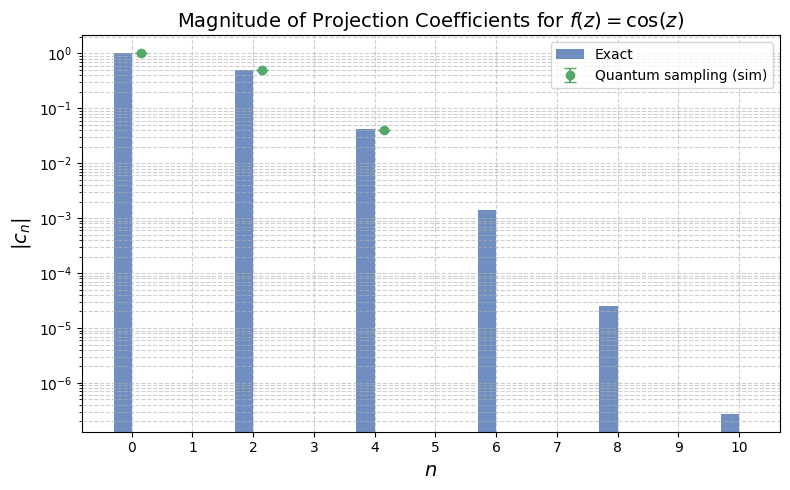}
\caption{Magnitude of the projection coefficients $|c_n|$ for $f(z) = \cos(z)$, comparing exact analytical values (blue bars) with simulated quantum measurement via Fock-state sampling (green circles) under realistic noise conditions ($\eta=0.95$, $10^5$ shots). The rapid factorial decay and strict parity (odd $n$ coefficients consistent with zero) are faithfully reproduced, validating the CV quantum algorithm for magnitude estimation.}
\label{fig:benchmark_cos}
\end{figure}

Figure~\ref{fig:benchmark_cos} compares simulated results against exact values for $\eta = 0.95$ (representative of transition-edge sensors \cite{Hadfield}). The magnitude estimate is obtained as $|c_n| = \sqrt{p_n}\,\|\ket{f}\|$, with statistical uncertainties derived from binomial sampling noise. 

\subsection{Derivation of Estimation Error Scaling}
\label{subsec:error_derivation}

To quantify the statistical uncertainty in reconstructing the state coefficients $\{c_n\}$ from experimental data, we consider the error propagation through the measurement and inversion process. Let the true probability of observing photon number $n$ be $p_n = |c_n|^2$. In an ideal detector with unity efficiency, the number of counts $k_n$ observed in $M$ shots follows a binomial distribution $B(M, p_n)$, which approximates a Poisson distribution for $p_n \ll 1$. The variance of the estimated probability $\hat{p}_n = k_n/M$ is given by standard shot noise:
\begin{equation}
\text{Var}(\hat{p}_n) = \frac{p_n(1-p_n)}{M} \approx \frac{p_n}{M},
\label{eq:ideal_variance}
\end{equation}
where the approximation holds for $p_n \ll 1$, typical for higher Fock states in weakly excited cavities.

Since the coefficient magnitude is $|c_n| = \sqrt{p_n}$, the error $\delta |c_n|$ is obtained via standard error propagation $\delta |c_n| = \left| \frac{d|c_n|}{dp_n} \right| \sqrt{\text{Var}(\hat{p}_n)}$:
\begin{equation}
\delta |c_n|_{\text{ideal}} = \frac{1}{2\sqrt{p_n}} \sqrt{\frac{p_n}{M}} = \frac{1}{2\sqrt{M}}.
\label{eq:ideal_error}
\end{equation}
This confirms the shot-noise-limited scaling $\propto 1/\sqrt{M}$ for ideal detection.

However, realistic detectors exhibit finite efficiency $\eta < 1$. The detection process is modeled as a binomial thinning channel where the measured photon number $m$ is related to the true number $n$ by the conditional probability
\begin{equation}
P(m|n) = \binom{n}{m} \eta^m (1-\eta)^{n-m}, \quad 0 \leq m \leq n.
\label{eq:loss_channel}
\end{equation}
Reconstructing the true probabilities $\{p_n\}$ from the measured frequencies $\{f_m\}$ requires inverting this loss matrix, which amplifies statistical fluctuations. For a state with mean photon number $\langle n \rangle = \sum_n n p_n$, the variance of the corrected count estimate acquires an excess noise term proportional to the loss rate $(1-\eta)$ . Specifically, to first order in $(1-\eta)$, the variance of the reconstructed probability scales as:
\begin{equation}
\text{Var}(\hat{p}_n)_{\eta} \approx \frac{p_n}{M} \left[ 1 + \frac{1-\eta}{\eta} \langle n \rangle \right]^2.
\label{eq:loss_variance}
\end{equation}
This result follows from propagating the binomial thinning statistics through the maximum-likelihood inversion of Eq.~\eqref{eq:loss_channel}. 

Applying the same propagation rule $\delta |c_n| = \frac{1}{2\sqrt{p_n}} \sqrt{\text{Var}(\hat{p}_n)_{\eta}}$, we obtain:
\begin{equation}
\delta |c_n| \approx \frac{1}{2\sqrt{p_n}} \sqrt{\frac{p_n}{M}} \left( 1 + \frac{1-\eta}{\eta} \langle n \rangle \right).
\end{equation}
Simplifying the $\sqrt{p_n}$ terms yields the final scaling relation:
\begin{equation}
\delta |c_n| \approx \frac{1}{2\sqrt{M p_n}} \left( 1 + \frac{1-\eta}{\eta} \langle n \rangle \right),
\label{eq:error_scaling}
\end{equation}
where $M$ is the number of measurement shots and $\langle n \rangle$ is the mean photon number. This expression confirms that the protocol maintains shot-noise-limited performance ($1/\sqrt{M}$ scaling), with detector inefficiency introducing a constant multiplicative overhead factor $\approx 1/\sqrt{\eta}$ for $\langle n \rangle \sim \mathcal{O}(1)$. For the parameters used in our simulation ($M = 10^5$ shots and $\eta = 0.95$), we achieve relative errors $\delta |c_n|/|c_n| \leq 1\%$ for the dominant coefficients ($n \leq 4$), which is adequate for benchmarking the protocol and analyzing low-order truncations.

The overhead factor $\left( 1 + \frac{1-\eta}{\eta} \langle n \rangle \right)$ becomes significant only for highly excited states with $\langle n \rangle \gg \eta/(1-\eta)$. For the bad-cavity QED benchmark considered here, $\langle n \rangle \lesssim 2$ and $\eta \geq 0.95$, so the efficiency-induced penalty remains below $10\%$. This validates the use of Eq.~\eqref{eq:error_scaling} for estimating statistical uncertainties in our numerical validation protocol.

The quantitative agreement is summarized in Table~\ref{tab:cos_coeffs}, where simulated magnitudes $|c_n|$ match exact values to within statistical error for dominant terms ($n = 0,2,4$), while odd-$n$ coefficients are consistently suppressed to zero—confirming faithful parity detection.

\begin{table}[htbp]
\centering
\caption{Magnitude of projection coefficients $|c_n|$ for $f(z) = \cos(z)$: exact analytical values versus simulated quantum measurement ($10^5$ shots) with realistic noise model ($\eta=0.95$, dark counts $\bar{n}_{\text{dark}}=10^{-3}$, preparation error $\sigma_{\text{prep}}=0.02$). Coefficients for odd $n$ vanish exactly due to even parity; small nonzero values in simulation are consistent with statistical noise.}
\label{tab:cos_coeffs}
\begin{tabular}{c c c}
\toprule
$n$ & $|c_n|$ (exact) & $|c_n|$ (quantum sim) \\
\midrule
0  & $1.0000 \times 10^{0}$  & $1.0003 \times 10^{0}$  \\
1  & $0.0000 \times 10^{0}$  & $0.0000 \times 10^{0}$  \\
2  & $5.0000 \times 10^{-1}$ & $4.9938 \times 10^{-1}$ \\
3  & $0.0000 \times 10^{0}$  & $0.0000 \times 10^{0}$  \\
4  & $4.1667 \times 10^{-2}$ & $4.0802 \times 10^{-2}$ \\
5  & $0.0000 \times 10^{0}$  & $0.0000 \times 10^{0}$  \\
6  & $1.3889 \times 10^{-3}$ & $0.0000 \times 10^{0}$  \\
7  & $0.0000 \times 10^{0}$  & $0.0000 \times 10^{0}$  \\
8  & $2.4802 \times 10^{-5}$ & $0.0000 \times 10^{0}$  \\
9  & $0.0000 \times 10^{0}$  & $0.0000 \times 10^{0}$  \\
10 & $2.7557 \times 10^{-7}$ & $0.0000 \times 10^{0}$  \\
\bottomrule
\end{tabular}
\end{table}

\section{Conclusion}
\label{sec:conclusion}

In this work, we have proposed and analyzed a continuous-variable quantum algorithm for estimating the projection coefficients of holomorphic functions within the Segal--Bargmann space. By leveraging the isometric isomorphism between the Segal--Bargmann space $\mathcal{H}_{\mathrm{SB}}$ and the Fock space of a single bosonic mode, we demonstrated that the projection coefficients $c_n$ coincide exactly with the Maclaurin series coefficients $a_n$ of the target function $f(z)$. This equivalence transforms the problem of series expansion into a quantum state tomography task, where coefficients are directly accessible as Fock-basis amplitudes.

We detailed a comprehensive protocol involving the construction of a state-preparation oracle $\hat{U}_f$, tailored to the analytic structure of the target function---ranging from deterministic photon-addition for polynomials to variational ansätze for general entire functions. The extraction of coefficients is achieved through photon-number-resolving detection for magnitudes and interferometric phase referencing for complex phases. Our numerical simulations on elementary functions, including $e^z$, $e^{iz}$, $\sin(z)$, and $\cos(z)$, validated the algorithm's accuracy, showing remarkable agreement with analytical Maclaurin expansions. Furthermore, the magnitude profiles of the coefficients were shown to effectively reveal intrinsic function properties, such as parity and leading-order terms, facilitating safe truncation for approximations.

Crucially, we analyzed the protocol's robustness under realistic experimental noise models, including detector inefficiency, dark counts, and state preparation infidelity. The derived error scaling indicates that the estimation remains shot-noise limited, with detector inefficiency introducing only a manageable multiplicative overhead for moderate mean photon numbers. Simulations incorporating these noise sources confirmed that high-fidelity coefficient estimation is feasible with current CV quantum hardware capabilities (e.g., $\eta \geq 0.95$).

This framework offers a hardware-native protocol for characterizing non-Gaussian states and analyzing functions defined by quantum oracles, providing a complementary approach to classical numerical integration. Future work may extend this methodology to multivariate functions via multi-mode CV systems and explore experimental realization on time-domain multiplexed optical platforms, paving the way for quantum-enhanced function analysis and spectral decomposition in complex physical systems.
Moreover, this work has implications for CV quantum implementation of probability theory through the relationship between probability theory and numerical analysis, which has been studied extensively over many decades since the work of Sul'din and Larkin \cite{Suldin, suldin1, Larkin}.

\section*{Acknowledgments}
The author would like to thank the reviewers for their insightful comments and constructive suggestions, particularly regarding the detailed analysis of oracle construction and the inclusion of realistic noise models, which have significantly strengthened this manuscript.

\vskip 5mm
{\bf Funding Declaration:}  Funding information is not applicable


\vskip 5mm
\begin{thebibliography}{99}
\bibitem{Ahlfors}
L. Ahlfors, \emph{Complex Analysis}, McGraw-Hill Education; 3rd edition (January 1, 1979)
\bibitem{Fokas}
M. J. Ablowitz and A. S. Fokas, \emph{Complex Variables: Introduction and Applications}, 2nd edition,
Cambridge University Press (2003).
\bibitem{Segal}
I. E. Segal, \emph{ Mathematical problems of relativistic physics}, in Kac, M. (ed.),
Proceedings of the Summer Seminar, Boulder, Colorado, 1960, Vol. II,
Lectures in Applied Mathematics, American Mathematical Society (1963).
\bibitem{Bargmann}
V. Bargmann, \emph{On a Hilbert space of analytic functions and an associated integral transform part I
},  Commun. Pure. Appl. Math. {\bf 14} (3): 187 (1961).
\bibitem{Hall}
B. C. Hall, \emph{Holomorphic Methods in Analysis and Mathematical Physics}, Contemporary Mathematics, Volume 260, pp. 1-59 (2000).
\bibitem{Perelomov}
A. Perelomov, \emph{Generalized Coherent States and Their Applications}, Springer Berlin, Heidelberg (1986).
\bibitem{Stenholm}
S. Stenholm, \emph{A Bargmann representation solution of the Jaynes—Cummings model},  Opt. Commun. {\bf 36}, Pages 75-78 (1981).
\bibitem{Bishop}
A. Vourdas, R. F. Bishop, \emph{Thermal coherent states in the Bargmann representation},  Phys. Rev. A {\bf 50}, 3331 (1994).
\bibitem{almasri}
M. W. AlMasri and  M. R. B. Wahiddin, \emph{Bargmann representation of quantum absorption refrigerators}, Rep. Math. Phys. {\bf 89} (2), Pages 185-198 (2022).
\bibitem{master}
M. W. AlMasri, M. R. B. Wahiddin, \emph{Quantum Decomposition Algorithm For Master Equations of Stochastic Processes: The Damped Spin Case
}, Mod. Phys. Lett. A, {\bf 37} (32), 2250216 (2022).
\bibitem{Wahiddin}
M. W. AlMasri, M. R. B. Wahiddin, \emph{Integral Transforms and PT-symmetric Hamiltonians}, Chinese Journal of Physics, Volume {\bf 85}, Pages 127-134 (2023)
\bibitem{neuron}
M. W. AlMasri, \emph{Multi-Valued Quantum Neurons}, Int.J. Theor.Phys {\bf 63}, 39 (2024).
\bibitem{complex}
M. W. AlMasri, \emph{On Logic Gates with Complex Numbers}, International Journal of Parallel, Emergent and Distributed Systems {\bf 39} (6), Pages 682-695 (2024).
\bibitem{central}
F. Shayegani, K.V. Samani, M. Abdi, M.W. AlMasri, \emph{ Dephasing dynamics of a central spin interacting with a classical Ising spin bath},  Physics Letters A {\bf 526}, 129966 (2024).
\bibitem{PRL}
S. D. Bartlett, B. C. Sanders, S. L. Braunstein, K. Nemoto, \emph{Efficient Classical Simulation of Continuous Variable Quantum Information Processes}, Phys. Rev. Lett., {\bf 88}, 097904 (2002).
\bibitem{review}
S. L. Braunstein and P. van Loock, \emph{Quantum information with continuous variables
}, Rev. Mod. Phys. {\bf 77}, 513 (2005).
\bibitem{Menicucci}
N. C. Menicucci, P. van Loock, M. Gu, C. Weedbrook, T. C. Ralph, and M. A. Nielsen, \emph{Universal Quantum Computation with Continuous-Variable Cluster States}, Phys. Rev. Lett. {\bf 97}, 110501 (2006).
\bibitem{Adesso}
G. Adesso and F. Illuminati, \emph{Entanglement in continuous-variable systems: recent advances and current perspectives}, J. Phys. A: Math. Theor. {\bf 40}, 7821 (2007).
\bibitem{Gu}
M. Gu, C. Weedbrook, N. C. Menicucci, T. C. Ralph, and P. van Loock, \emph{Quantum computing with continuous-variable clusters}, Phys. Rev. A {\bf 79}, 062318 (2009).
\bibitem{scipy}
P. Virtanen \textit{et al.}, SciPy 1.0 Contributors, \textit{Nat. Methods} \textbf{17}, 261 (2020).
\bibitem{killoran}
N. Killoran,  T. R. Bromley, J. M. Arrazola, M. Schuld, N. Quesada, and S. Lloyd \emph{Continuous-variable quantum neural networks}, Phys. Rev. Res. {\bf 1}, 033063 (2019).
\bibitem{Zavatta2004}
A. Zavatta, S. Viciani, and M. Bellini, \emph{Tomographic reconstruction of the single-photon Fock state by high-frequency homodyne detection}, Phys. Rev. A {\bf 70}, 053821 (2004).
\bibitem{nelder1965simplex} 
J.A. Nelder, R. Mead, \emph{A simplex method for function minimization}, The Computer Journal, 1965, vol. 7, no. 4, 308--313
\bibitem{Menicucci2014}
N. C. Menicucci, \emph{Fault-tolerant measurement-based quantum computing with continuous-variable cluster states}, Phys. Rev. Lett. {\bf 112}, 120504 (2014).
\bibitem{tomography}
A. I. Lvovsky, H. Hansen, T. Aichele, O. Benson, J. Mlynek, S. Schiller, \emph{Quantum state reconstruction of the single-photon Fock state}, Phys. Rev. Lett. {\bf 87}, 050402 (2001).
\bibitem{tomography1}
A. Zavatta, S. Viciani and M. Bellini, \emph{Tomographic reconstruction of the single-photon Fock state by high-frequency homodyne detection}, Phys. Rev. A {\bf 70}, 053821 (2004).



\bibitem{Hadfield}
R. Hadfield, \emph{ Single-photon detectors for optical quantum information applications}. Nature Photon {\bf 3}, 696–705 (2009).
\bibitem{Suldin}
A.V. Sul'din, \emph{ Wiener measure and its applications to approximation methods. I.}  Izv. Vyssh. Uchebn. Zaved. Mat., 1959, no. 6, 145–158
\bibitem{suldin1}
A.V. Sul'din, \emph{ Wiener measure and its applications to approximation methods. II.},  Izv. Vyssh. Uchebn. Zaved. Mat., 1960, no. 5, 165–179
\bibitem{Larkin}
F.M. Larkin, \emph{Gaussian measure in Hilbert space and applications in numerical analysis}, Rocky Mountain J. Math. 2(3): 379-422 (1972).
\end{thebibliography}
\end{document}